\documentclass[reqno,hidelinks,a4paper]{amsart}
\usepackage[T1]{fontenc}

\usepackage{latexsym, amssymb, amsmath}
\usepackage{enumitem}
\usepackage{amsthm}
\usepackage{tikz-cd} 
\usepackage{amsthm} 
\usepackage{orcidlink}

\usepackage[utf8]{inputenc}
\usepackage[british]{babel}
\usepackage[all]{xy}
\usepackage{hyperref}
\usepackage{calrsfs}
\usepackage{mathabx}
\usepackage{xcolor}
\usepackage{mathtools}
\usepackage{blkarray}
\usepackage[textsize=scriptsize]{todonotes}
\usepackage{thmtools,mathtools}
\usepackage{nameref,hyperref,cleveref}

\theoremstyle{plain}
\newtheorem{theorem}[subsection]{Theorem}
\newtheorem{proposition}[subsection]{Proposition}
\newtheorem{corollary}[subsection]{Corollary}

\theoremstyle{definition}
\newtheorem{definition}[subsection]{Definition}

\theoremstyle{remark}

\numberwithin{equation}{section}

\newcommand{\noproof}{\hfil\qed}

\newcommand{\C}{\ensuremath{\mathcal{C}}}

\newcommand{\F}{\mathbb{F}}

\newcommand{\V}{\mathcal{V}}
\newcommand{\Z}{\mathbb{Z}}

\DeclareMathOperator{\Aut}{Aut}

\DeclareMathOperator{\Act}{Act}
\DeclareMathOperator{\Hom}{Hom}
\DeclareMathOperator{\Bim}{Bim}
\DeclareMathOperator{\SplExt}{SplExt}
\DeclareMathOperator{\Imm}{Im}

\DeclareMathOperator{\Der}{Der}

\DeclareMathOperator{\Bider}{Bider}

\DeclareMathOperator{\id}{id}

\newcommand{\AAlg}{\ensuremath{\mathbf{Assoc}}}

\newcommand{\Lie}{\ensuremath{\mathbf{Lie}}}
\newcommand{\LeibAlg}{\ensuremath{\mathbf{Leib}}}

\newcommand{\CAAlg}{\ensuremath{\mathbf{CAssoc}}}

\newcommand{\AbAlg}{\ensuremath{\mathbf{AbAlg}}}

\newcommand{\Nil}{\ensuremath{\mathbf{Nil}}}
\newcommand{\Sol}{\ensuremath{\mathbf{Sol}}}

\newcommand{\Grp}{\ensuremath{\mathbf{Grp}}}

\newdir{>}{{}*:(1,-.2)@^{>}*:(1,+.2)@_{>}}
\newdir{<}{{}*:(1,+.2)@^{<}*:(1,-.2)@_{<}}
\newdir{>>}{{}*!/3.5pt/:(1,-.2)@^{>}*!/3.5pt/:(1,+.2)@_{>}*!/7pt/:(1,-.2)@^{>}*!/7pt/:(1,+.2)@_{>}}
\newdir{ >}{{}*!/-8pt/@{>}}

\date{}

\begin{document}
	
\title[Action accessible and weakly action repr.~varieties of algebras]{Action accessible and weakly action representable varieties of algebras}
	
\author[X.~García-Martínez]{X.~García-Martínez~\orcidlink{0000-0003-1679-4047}}
\author[M.~Mancini]{M.~Mancini$^{*}$~\orcidlink{0000-0003-2142-6193}}
	
\email{xabier.garcia@usc.es}
\email{manuel.mancini@unipa.it; manuel.mancini@uclouvain.be}
	
\address[Xabier García-Martínez]{CITMAga \& Facultade de Matemáticas, Universidade de Santiago de Compostela, Rúa de Lope Gómez de Marzoa, s/n, 15782, Santiago de Compostela, Spain.}
\address[Manuel Mancini]{Dipartimento di Matematica e Informatica, Universit\`a degli Studi di Palermo, Via Archirafi 34, 90123 Palermo, Italy.}
\address[Manuel Mancini]{Institut de Recherche en Mathématique et Physique, Université catholique de Louvain, chemin du cyclotron 2 bte L7.01.02, B--1348 Louvain-la-Neuve, Belgium.}

\thanks{The first author is supported by Ministerio de Economía y Competitividad (Spain) with grant number PID2021-127075NA-I00. The second author is supported by the University of Palermo; by the ``National Group for Algebraic and Geometric Structures and their Applications'' (GNSAGA -- INdAM); by the National Recovery and Resilience Plan (NRRP), Mission 4, Component 2, Investment 1.1, Call for tender No.~1409 published on 14/09/2022 by the Italian Ministry of University and Research (MUR), funded by the European Union -- NextGenerationEU -- Project Title Quantum Models for Logic, Computation and Natural Processes (QM4NP) -- CUP:~B53D23030160001 -- Grant Assignment Decree No.~1371 adopted on 01/09/2023 by the Italian Ministry of Ministry of University and Research (MUR); by the SDF Sustainability Decision Framework Research Project -- MISE decree of 31/12/2021 (MIMIT Dipartimento per le politiche per le imprese -- Direzione generale per gli incentivi alle imprese) -- CUP:~B79J23000530005, COR:~14019279, Lead Partner:~TD Group Italia Srl, Partner:~University of Palermo; and he is also a postdoctoral researcher of the Fonds de la Recherche Scientifique--FNRS}
	
\begin{abstract}
The main goal of this article is to investigate the relationship between action accessibility and weak action representability in the context of varieties of non-associative algebras over a field. Specifically, using an argument of J.~R.~A.~Gray in the setting of groups, we prove that the varieties of $k$-nilpotent Lie algebras ($k \geq 3$) and the varieties of $n$-solvable Lie algebras ($n \geq 2$) do not form weakly action representable categories. These are the first known examples of action accessible varieties of non-associative algebras that fail to be weakly action representable, establishing that a subvariety of a (weakly) action representable variety of non-associative algebras needs not be weakly action representable. Eventually, we refine J.~R.~A.~Gray's result by proving that the varieties of $k$-nilpotent groups ($k \geq 3$) and that of $2$-solvable groups are not weakly action representable.
\end{abstract}
	
	
\maketitle
	
\section*{Introduction}\label{Introduction}
	
\emph{Internal actions} were defined in \cite{IntAct} by F.~Borceux, G.~Janelidze and G.~M.~Kelly with the aim of extending the correspondence between actions and split extensions from the context of groups to arbitrary semi-abelian categories~\cite{Semi-Ab}. In some of those categories, internal actions are exceptionally well behaved, in the sense that the actions on each object~$X$ are \emph{representable}: this means that there exists an object~$[X]$ such that the functor $\Act(-,X)\cong \SplExt(-,X)$ which sends an object $B$ to the set of actions of $B$ on $X$ (isomorphisms classes of split extensions of $B$ by~$X$) is naturally isomorphic to the functor $\Hom(-,[X])$. The notion of representability of actions in a semi-abelian category is further studied in~\cite{BJK2}, where it is explained for instance that the category of commutative associative algebras over a field is not action representable. Later it was shown in~\cite{Tim} that the only action representable varieties of non-associative algebras over an infinite field of characteristic different from $2$ are the category $\Lie$ of Lie algebras and the category $\AbAlg$ of abelian algebras. The relative strength of the notion naturally led to the definition of closely related weaker notions.
	
In~\cite{act_accessible} D.~Bourn and G.~Janelidze introduced the notion of \emph{action accessible} category in order to include relevant examples that do not fit into the frame of action representable categories (such as rings, associative algebras and Leibniz algebras~\cite{loday}). A.~Montoli proved in~\cite{Montoli} that any \emph{category of interest} in the sense of Orzech~\cite{Orzech} is action accessible, while in~\cite{Casas} the authors showed that a weaker notion of representing object (which they name as the \emph{universal strict general actor}) is available for any Orzech category of interest. 
	
Recently, G.~Janelidze introduced in~\cite{WRA} the concept of \emph{weakly action representable} category. Instead of asking that for any object $X$ in a semi-abelian category~$\C$ there is an object $[X]$ and a natural isomorphism $\Act(-,X)\cong\Hom_\C(-,[X])$, we require the existence of a \emph{weakly representing object} $T=T(X)$ and a natural monomorphism of functors
\[
\tau\colon\Act(-,X)\rightarrowtail\Hom_\C(-,T).
\]
When each $X$ admits such an object $T$, $\C$ is said to be \emph{weakly action representable}. This is the case for instance of the category of associative algebras~\cite{WRA}, the category of commutative associative algebras~\cite{WRAAlg}, or the category of Leibniz algebras~\cite{CigoliManciniMetere}.
	
In~\cite{WRA} the author proved that any weakly action representable category is action accessibile. However, J.~R.~A.~Gray observed in~\cite{Gray} that the converse of this implication is not true: he showed that the notion of weak action representability is connected to the existence of a so-called \emph{amalgam}~\cite{kiss}, which already appeared in~\cite{BJK2} in relations with action representability, and he proved that the varieties of $n$-solvable groups ($n \geq 3$) are not weakly action representable.
	
The notion of weakly representable action was then studied in \cite{WRAAlg} in the context of varieties of non-associative algebras over a field, where the authors provided new examples of weakly action representable categories, such as the variety of commutative associative algebras and those of $2$-nilpotent (commutative, anti-commutative and non-commutative) algebras. 
	
Nevertheless, the authors were not able to find an example of an action accessible variety which is not weakly action representable and they ended up with the following open questions: 
\begin{enumerate}
\item Does the converse of the implication
\[
\textit{weakly action representable category} \Rightarrow \textit{action accessible category}
\]
hold in the context of varieties of non-associative algebras over a field?
\item How does the condition of weakly representable actions behave under taking subvarieties? Specifically, if a variety of non-associative algebras $\V$ is weakly action representable, does it follow that every subvariety of $\V$ is also weakly action representable?
\end{enumerate}
	
The main aim of this article is to show that the answer to these questions is ``no''. We use Proposition 2.2.~of~\cite{Gray} to prove that the varieties $\Nil_k(\Lie)$ of $k$-nilpotent Lie algebras and $\Sol_n(\Lie)$ of $n$-solvable Lie algebras, which are subvarieties of $\Lie$, are not weakly action representable for any $k \geq 3$ and $n \geq 2$. As an additional result, we use similar methods to refine J.~R.~A.~Gray's result~\cite{Gray} by proving that the varieties of $k$-nilpotent groups ($k \geq 3$) and that of $2$-solvable groups are not weakly action representable.
	
\section{Preliminaries}\label{SecPrel}
	
Let $\C$ be a semi-abelian category~\cite{Semi-Ab} and let $B,X$ be objects of $\C$. We recall that a \emph{split extension} of $B$ by $X$ is a diagram
	
\begin{equation*}\label{eq:split_ext}
\begin{tikzcd}
0\ar[r]
&X\arrow [r, "k"]
&A \arrow[r, shift left, "\alpha"] & 
B\ar[r]\ar[l, shift left, "\beta"]
&0 
\end{tikzcd}	
\end{equation*}
in $\C$ such that $\alpha \circ \beta = \id_B$ and $(X,k)$ is a kernel of $\alpha$.

For any object $X$ of $\C$, one may define the functor
\[
\SplExt(-,X)\colon \C^{op} \rightarrow \textbf{Set}
\]
which maps any object $B$ of $\C$ to the set $\SplExt(B,X)$ of isomorphism classes of split extensions of $B$ by $X$, and any arrow $f\colon B'\to B$ to the \emph{change of base} morphism $f^*\colon \SplExt(B,X) \to \SplExt(B',X)$ given by pulling back along $f$.
	
Given a semi-abelian category, one may define the notion of internal action~\cite{IntAct}. Internal actions on an object $X$ give rise to a functor 
\[
\Act(-,X)\colon \C^{op} \rightarrow \textbf{Set}
\]
and one may prove there exists a natural isomorphism $\Act(-,X)\cong \SplExt(-,X)$ (see~\cite{BJK2}). Here we don't describe explicitly internal actions, since split extensions are more handy to work with, especially in the context of non-associative algebras. This justifies the terminology in the definition that follows. 
	
\begin{definition}[\cite{BJK2}]
A semi-abelian category $\C$ is said to be \emph{action representable} if for any object $X$ in $\C$, the functor $\SplExt(-,X)$ is representable. This means that there exists an object $[X]$ of $\C$ and a natural isomorphism of functors
\[
\SplExt(-,X) \cong \Hom_{\C}(-,[X]).
\]
\end{definition}
	
The prototype examples of action representable categories are the category $\Grp$ of groups and the category $\Lie$ of  Lie algebras over a commutative unitary ring. In the first case, it is well known that every action of $B$ by $X$ is represented by a group homomorphism $B \rightarrow \Aut(X)$, where $\Aut(X)$ is the group of automorphisms of $X$. In the second case, any split extension of $B$ by $X$ is represented by a Lie algebra homomorphism $B \rightarrow \Der(X)$, where $\Der(X)$ is the Lie algebra of derivations of $X$.
	
However the notion of action representable category has proven to be quite restrictive: for instance, in~\cite{Tim} the authors proved that the only examples of action representable varieties of non-associative algebras (over an infinite field $\F$ of characteristic different from $2$) are the category $\AbAlg$ of abelian algebras and the category $\Lie$ of Lie algebras.
		
\begin{definition}[\cite{WRA}]
A semi-abelian category $\C$ is said to be \emph{weakly action representable} if for every objext $X$ in $\C$, the functor $\SplExt(-,X)$ admits a weak representation. This means that there exist an object $T=T(X)$ of $\C$ and a natural monomorphism of functors
\[
\tau\colon \SplExt(-,X) \rightarrowtail \Hom_{\C}(-,T).
\]
A morphism $(\varphi\colon B \rightarrow T) \in \Imm(\tau_B)$ is called \emph{acting morphism}.
\end{definition}
	
Examples of weakly action representable category include the variety $\AAlg$ of associative algebras~\cite{WRA}, where $T=\Bim(X)$ is the associative algebra of \emph{bimultipliers} of $X$~\cite{MacLane58}; the variety $\LeibAlg$ of Leibniz algebras, where $T=\Bider(X)$ is the Leibniz algebra of \emph{biderivations} of $X$~\cite{loday,ManciniBider}; the varieties of $2$-nilpotent (commutative, anti-commutative and non-commutative) algebras (see \cite{WRAAlg, LaRosaMancini1, LaRosaMancini2, LaRosaMancini3}). Note that in these mentioned examples, the weakly representing object is quite easy to consider. However, this is not always the case. The variety $\CAAlg$ of commutative associative algebras is also weakly action representable, but the algebra that could be expected to be a weakly representing object is not commutative. Nevertheless, in~\cite{WRAAlg} a big colimit merged with the amalgamation property is used to naturally construct a weakly representing object.
	
Another important observation made by G.~Janelidze in~\cite{WRA} is that every weakly action representable category is action accessible. We thus have that
\[
\textit{action representability}\Rightarrow\textit{weak action representability}\Rightarrow\textit{action accessibility.}
\]
The author ended up with an open question, whether reasonably mild conditions may be found on a semi-abelian category under which the second implication may be reversed.
	
This question was very interestingly addressed by J.~R.~A.\ Gray, who proved in~\cite{Gray} that the varieties of $n$-solvable groups are not weakly action representable for any $n \geq 3$. He showed that the notion of weak action representability is connected with the so-called \emph{amalgamation property} (AP)~\cite{kiss}.
	
We recall that a span of monomorphisms $m \colon S \rightarrowtail B$ and $m' \colon S \rightarrowtail B'$ in a category $\C$ admits an amalgam in $\C$ if there exist two monomorphisms $u \colon B \rightarrowtail D$ and $u' \colon B' \rightarrowtail D$ in $\C$ such that $um=u'm'$. The relations between the representabilty of actions and the (AP) were firstly investigated in~\cite{BJK2}, where it was proved that for any Orzech category of interest, action representability is equivalent to the (AP) for \emph{protosplit} monomorphisms, i.e., monomorphisms which form the kernel part of a split extension.
	
Using the (AP), J.~R.~A.\ Gray gave sufficient conditions under which a Birkhoff subcategory of an action representable category is not weakly action representable.
	
\begin{proposition}[\cite{Gray}]\label{Gray}
Let $\C$ be an action representable category and let $\mathcal{X}$ be a Birkhoff subcategory of $\C$. Suppose there exist two monomorphisms $m \colon S \rightarrowtail B$, $m' \colon S \rightarrowtail B'$ in $\mathbb{X}$, two monomorphisms $u \colon B \rightarrowtail D$, $u' \colon B' \rightarrowtail D$ in $\C$, an object $X$ of $\mathcal{X}$ and a monomorphism $v \colon D \rightarrowtail [X]$ such that
\begin{itemize}
\item[(a)] $m$ and $m'$ cannot be amalgamated in $\mathcal{X}$;
\item[(b)] $um=u'm'$;
\item[(c)] the split extensions with kernel $X$ corresponding to $vu$ and $vu'$ are in $\mathcal{X}$. 
\end{itemize}
Then, the category $\mathcal{X}$ is not weakly action representable. \noproof
\end{proposition}
	
Thanks to a B.~H.~Neumann's example of an abelian group $S$, a $2$-nilpotent group $B$ and two monomorphisms $m,m' \colon S \rightarrowtail B$ which cannot be amalgamated in any solvable group $D$ (see~\cite{Neumann}), it is possible to apply \Cref{Gray} in order to prove that the varieties of $n$-solvable groups are not weakly action representable for any $n \geq 3$.



\section{Nilpotent and \texorpdfstring{$2$}{2}-solvable groups}

In this section we adapt to $k$-nilpotent ($k \geq 3)$ and $2$-solvable groups the proof stating that the categories of $n$-solvable groups ($n \geq 3$) are not weakly action representable. 

Let us first recall the concrete example the example of B.~H.~Neumann that is used in J.~R.~A.~Gray's proof. Consider the following groups in the form of generators and relations:
\begin{align*}
S &= \langle a, b \; | \; [a, b]=a^5=b^5=1 \rangle , \\ 
B &= \langle x, a, b \; | \; [x, a]=b^{-1}, [x, b] = x^5 = a^5 = b^5=1 \rangle, \\
B'&= \langle y, a, b \; | \; [y, b] = a^{-1}, [y, a] = y^5 = a^5 = b^5 = 1\rangle.
\end{align*}

It was proved in~\cite{Neumann} that the obvious inclusions of $S$ inside $B$ and $B'$ cannot be amalgamated in any solvable group, and it was observed in~\cite{Gray} that there exist a group $X$ and two monomorphisms $\psi \colon B \rightarrowtail \Aut(X)$ and $\psi' \colon B' \rightarrowtail \Aut(X)$ agreeing on $S$, such that the corresponding split extensions are $3$-solvable.

We aim now to show that it is possible to choose $X$, $\psi$ and $\psi'$ in such a way the split extensions corresponding to $\psi$ and $\psi'$ are $3$-nilpotent. Let $X=\Z_5^3$ be the abelian group formed by $3$ copies of $\Z_5$. We consider the group monomorphisms $\psi \colon B \rightarrowtail \Aut(X)$ and $\psi' \colon B' \rightarrowtail \Aut(X)$ defined by
\begin{align*}
\psi(x) &= 
\begin{pmatrix}
1 & 1 & 0\\
0 & 1 & 0\\
0 & 0 & 1
\end{pmatrix}, &
\psi(y) &= 
\begin{pmatrix}
1 & 0 & 0\\
1 & 1 & 0\\
0 & 0 & 1
\end{pmatrix},
\\ \\
\psi(a) = \psi'(a) &=
\begin{pmatrix}
1 & 0 & 0\\
0 & 1 & 1\\
0 & 0 & 1
\end{pmatrix},&
\psi(b) = \psi'(b) &= 
\begin{pmatrix}
1 & 0 & 4\\
0 & 1 & 0\\
0 & 0 & 1
\end{pmatrix}.
\end{align*}

It is immediate to see that the relations are preserved. To finish the proof, we just need to check that the induced split extensions are 3-nilpotent. 

Let $B_\psi=B \ltimes_\psi X$ and $B_{\psi'}=B \ltimes_{\psi'} X$ be the semi-direct products induced by $\psi$ and $\psi'$ respectively. We have to check that they are $3$-nilpotent groups. Recall that the multiplication in $B_\psi$ is defined by 
\[
(g,x) (g',x')=(gg',x \psi(g)(x')),
\]
for any $g,g' \in B$ and $x,x' \in X$. Since $B$ is $2$-nilpotent and $X$ is abelian, one may check that
\begin{align*}
&[B_\psi,B_\psi] \cong\langle b \rangle \ltimes_\psi \mathbb{Z}_5^2,\\
&[B_\psi,[B_\psi,B_\psi]] \cong \mathbb{Z}_5    
\end{align*}
and $[B_\psi,[B_\psi,[B_\psi,B_\psi]]]$ is the trivial group.
Thus,~$B_\psi$ is a $3$-nilpotent group. The same result can be obtained for the semi-direct product $B_{\psi'}$ and we can state the following.

\begin{theorem}
The variety of $k$-nilpotent groups is not weakly action representable for any $k \geq 3$. \noproof
\end{theorem}

Since $3$-nilpotent implies $2$-solvable, the same example allows to prove the following.

\begin{corollary}
The variety of $2$-solvable groups is not weakly action representable. \noproof
\end{corollary}
    
\section{Subvarieties of Lie algebras}\label{Solvable}
The aim of this section is to prove that, over a field $\F$, the varieties of $k$-nilpotent Lie algebras ($k \geq 3$) and the varieties of $n$-solvable Lie algebras ($n \geq 2$) are not weakly action representable. 
	
Recall that, given any Lie algebra $L$, its lower central series is defined as
\[
L^{0}=L, \quad L^{k}=[L,L^{k-1}],
\]
for any $k \in \mathbb{N}$. The Lie algebra $L$ is said to be $k$-nilpotent if $L^{k-1} \neq 0$ and $L^{k}=0$. In addition, its derived series is defined as
\[
L^{(0)}=L, \quad L^{(n)}=[L^{(n-1)},L^{(n-1)}],
\]
for any $n \in \mathbb{N}$. The Lie algebra $L$ is said to be $n$-solvable if $L^{(n-1)} \neq 0$ and $L^{(n)}=0$.
	
We denote by \( \Nil_k(\Lie) \) and \( \Sol_n(\Lie) \) the subvarieties of \( \Lie \) consisting of all Lie algebras that are \( s \)-nilpotent and \( t \)-solvable, respectively, for some \( s \leq k \) and \( t \leq n \). Note that the classes of all nilpotent and solvable Lie algebras do not, in general, form varieties; however, they do when the degree of solvability or nilpotency is bounded.
	
Inspired by B.~H.~Neumann’s proof for solvable groups \cite{Neumann}, we show that there exist an abelian algebra $S$, and two monomorphisms $m\colon S \rightarrowtail B$ and $m'\colon S \rightarrowtail B'$, where $B$ and $B'$ are solvable, which cannot be amalgamated in any solvable Lie algebra.
	
Let $S=\F \lbrace a,b \rbrace$ be the $2$-dimensional abelian algebra and consider two copies $B, B'$ of the $3$-dimensional Heisenberg algebra. More explicitly, $B= \F \lbrace x, a, b\rbrace$ with the only non-zero bracket $[x,a]=b$, and $B'= \F \lbrace y, a, b\rbrace$ with the only non-zero bracket $[y,b]=a$, where the mononomorphisms are the obvious inclusions.
	
Now, suppose there exists a solvable Lie algebra $D$ and two monomorphisms $u \colon B \rightarrowtail D$ and $u'\colon B' \rightarrowtail D$, such that they agree on $S$.
	
Consider the Lie subalgebra $P$ of $D$ generated by $u(B)$ and $u'(B')$. Since $S$ is an ideal of $B$ and $B'$, then $U=u(S)=u'(S)$ is an ideal of $P$. This gives us the adjoint map
\[
\operatorname{ad} \colon P \rightarrow \Der(U), \quad p \mapsto \operatorname{ad}_p,
\]
where $\operatorname{ad}_p=[p,-]$ and $\Der(U) \cong \mathfrak{gl}(2,\F)$. 
	
Since $P$ is solvable, $\operatorname{ad}(P)\cong P/\ker(\operatorname{ad})$ is also solvable. Nevertheless, $\operatorname{ad}(P)$ is the Lie algebra generated by 
\[
\operatorname{ad}_{u(x)} \equiv \begin{pmatrix}
0 & 0 \\
1 & 0 
\end{pmatrix} \text{ and } \operatorname{ad}_{u'(y)} \equiv
\begin{pmatrix}
0 & 1 \\
0 & 0
\end{pmatrix}
\]

Hence, $\operatorname{ad}(P)$ is isomorphic to the special linear algebra $\mathfrak{sl}(2,\F)$, which is a simple Lie algebra, and we get a contradiction. Thus, $m$ and $m'$ cannot be amalgamated in any solvable Lie algebra. Since nilpotency implies solvability, this argument shows that the monomorphisms $m \colon S \rightarrowtail B$ and $m' \colon S \rightarrowtail B'$ cannot be amalgamated in any nilpotent Lie algebra. To conclude the proof, we observe that
\begin{itemize}
\item[(i)] $\Lie$ has the amalgamation property~\cite{kiss}. Thus, there exist a Lie algebra $L$ and two monomorphisms $\psi \colon B \rightarrowtail L$ and $\psi' \colon B' \rightarrowtail L$ agreeing on $S$. For instance, we may consider $L=\mathfrak{gl}(3,\F)$ and the two faithful representations defined by
\[
\psi(x)=-e_{23}, \quad \psi(a)=e_{12}, \quad \psi(b)=e_{13}
\]
and
\[
\psi'(y)=-e_{32}, \quad , \psi'(a)=e_{12} \quad, \psi'(b)=e_{13},
\]
where $e_{ij}$ is the $3 \times 3$ matrix which has $1$ in the entry $(i,j)$ and $0$ elsewhere. It is immediate to check that the subalgebras $\psi(B)$ and $\psi'(B')$ of $\mathfrak{gl}(3,\F)$ are both isomorphic to the $3$-dimensional Heisenberg algebra.
\item[(ii)] Let $X$ be the $3$-dimensional abelian algebra and let $v=\id_{L}$ be the canonical isomorphism $L \cong \Der(X)$. We aim to prove that the split extension of $B$ by $X$ corresponding to $v\circ \psi=\psi$ and that of $B'$ by $X$ corresponding to $v \circ \psi' = \psi'$ are in $\Nil_k(\Lie)$, for any $k \geq 3$. To get the result, we need to explicitly check that the semi-direct products $B_\psi= B \ltimes_\psi X$ and $B'_{\psi'}= B \ltimes_{\psi'} X$, induced by $\psi$ and $\psi'$ respectively, are $3$-nilpotent Lie algebras. 
In fact, the semi-direct product $B_\psi$ may be described as the Lie algebra with basis $\lbrace x,a,b,e_1,e_2,e_3 \rbrace$ and Lie brackets
\[
[x,a]=b, \quad [x,e_3]=-e_2, \quad [a,e_2]=[b,e_3]=e_1.
\]
One may easily check that $[B_\psi,B_\psi]=\F \lbrace b,e_1,e_2 \rbrace$, $[B_\psi,[B_\psi,B_\psi]]=\F \lbrace e_1 \rbrace$ and $[B_\psi,[B_\psi,[B_\psi,B_\psi]]]=0$, which means that $B_\psi$ is $3$-nilpotent. This can be done similarly in the case of $B'_{\psi'}$.
\end{itemize}
	
Summarising, the hypothesis of \Cref{Gray} are verified and we can conclude with the following.

\begin{theorem}
The variety $\Nil_k(\Lie)$ is not weakly action representable for any $k \geq 3$. \noproof
\end{theorem}

Note that the same argument does not apply to the case $k=2$ since the semi-direct products $B_\psi$ and $B_{\psi'}$ are not $2$-nilpotent algebras. 

Finally, since any $3$-nilpotent Lie algebras is also $2$-solvable, we get the following.

\begin{corollary}
The variety $\Sol_n(\Lie)$ is not weakly action representable for any $n \geq 2$. \noproof
\end{corollary}

\section*{Acknowledgements}
The authors would like to thank the anonymous referees for their useful comments that helped us to improve the quality of the manuscript.


\end{document}